\documentclass[11pt,a4paper]{article}
\usepackage[T1]{fontenc}
\usepackage[utf8]{inputenc}
\usepackage[margin=1in]{geometry}
\usepackage{lmodern}
\usepackage{amsmath,amssymb,amsfonts,amsthm}
\usepackage{graphicx,subcaption,mathrsfs,dsfont,bm,enumerate,xcolor}
\usepackage{array,multirow,rotating,float}
\usepackage[section]{placeins}
\usepackage[authoryear,round]{natbib}
\usepackage{microtype}
\usepackage{setspace}
\usepackage{hyperref}
\hypersetup{
  colorlinks=true,
  linkcolor=blue,
  citecolor=blue,
  urlcolor=blue,
  pdftitle={Testing Equality of Distributions via Repeatedly Integrated Quantile Metrics Under Weak Moment Conditions},
  pdfauthor={Zhenfeng Zou, Meng Guan, Panxu Yuan}
}

\newtheorem{theorem}{Theorem}[section]

\newtheorem{corollary}[theorem]{Corollary}
\newtheorem{proposition}[theorem]{Proposition}
\theoremstyle{definition}
\newtheorem{definition}{Definition}[section]

\newtheorem{condition}{Condition}

\numberwithin{equation}{section}

\def\oo{\infty} \def\d{\,\mathrm{d}}  
\newcommand{\dist}{\mathrm{d}}
 
\newcommand{\esssup}{\mathrm{ess\mbox{-}sup}}

\newcommand{\E}{\mathsf{E}}
\newcommand{\Prob}{\mathsf{P}}
\newcommand{\Var}{\mathsf{Var}}
\newcommand{\eqd}{\overset{\dist}{=}}
\newcommand{\dto}{\xrightarrow{\dist}}
\newcommand{\pto}{\xrightarrow{\Prob}}
\newcommand{\asto}{\xrightarrow{\mathrm{a.s.}}}
\newcommand{\sgn}{\operatorname{sgn}}

\newcommand{\R}{\mathbb{R}}

\renewcommand{\(}{\left(}
\renewcommand{\)}{\right)}
\renewcommand{\[}{\left[}

\newcommand{\1}{\mathds{1}}
\allowdisplaybreaks

\title{ Testing Equality of Distributions via Repeatedly Integrated Quantile Metrics Under Weak Moment Conditions}
\author{
Zhenfeng Zou\thanks{School of Public Policy and Management, University of Science and Technology of China,  China. Email: \url{zfzou@ustc.edu.cn}}
\and Meng Guan\thanks{School of Management, University of Science and Technology of China, China.
Email: \url{gm123123@mail.ustc.edu.cn}}
\and Panxu Yuan\thanks{Corresponding author, School of Management, University of Science and Technology of China,  China. Email: \url{yuanpanxu@foxmail.com}}
\and Sanying Feng\thanks{Co-corresponding author, School of Mathematics and Statistics, Zhengzhou University,  China. Email: \url{fsy5801@zzu.edu.cn}}}
\date{}

\begin{document}
\maketitle

\begin{abstract}
Testing whether two independent samples arise from the same underlying distribution is a fundamental statistical problem. We propose a new class of two-sample distribution tests based on a family of probability metrics $\Delta_{n,p}$, constructed from repeatedly integrated quantile functions. On their respective domains, these metrics are proved to be genuine distributional distances. The case $n=1$ recovers the $p$-Wasserstein distance, which requires finite $p$-th moments; for $n\geq2$,  the proposed metrics are well defined and require only finite first moments. The asymptotic properties of the plug-in statistic are established, including strong consistency and limiting distributions under the null and fixed alternatives. A permutation calibration for finite-sample inference is also proposed. We further derive an asymptotic power function under local alternatives. Finally, the finite-sample performance of the proposed tests is examined through simulation studies, and their reduced sensitivity to extreme upper-tail observations is illustrated through a real data application.
\end{abstract}


\section{Introduction}

Testing the equality of two distributions from independent samples is a basic problem in nonparametric statistics \citep{LehmannRomano22}. Such comparisons are useful when two populations may differ in variability or tail behavior, even if a comparison of their means is not informative. Examples include treatment and control outcomes in clinical trials \citep{Pocock13,Fried15}, earnings distributions in program evaluation \citep{Bloom97,AAI02}, and heavy-tailed financial returns or insurance losses \citep{Cont01,Embrechts97}. These settings motivate two-sample procedures that compare the distributions as a whole under weak moment assumptions.

The two-sample problem has a long history, and many nonparametric procedures are available,
for example, the Kolmogorov--Smirnov test \citep{Kolmogorov33,Smirnov39},
the Cram\'{e}r--von Mises test \citep{Anderson62},
the Anderson--Darling test \citep{Pettitt76},
and the Wilcoxon--Mann--Whitney test \citep{Wilcoxon45,Mann47}.
More recent procedures include the energy-distance test \citep{Baringhaus04,SzekelyRizzo13} and the maximum mean discrepancy (MMD) test \citep{Gretton12}. These procedures summarize distributional differences in different ways and therefore respond differently to different alternatives. For a broader discussion of these and related two-sample procedures, see \citet{RGC17} and the references therein.

Among quantile-based approaches for distributions on the real line, the Wasserstein distance provides a natural comparison through corresponding quantile levels \citep{Vallender74}. Its empirical asymptotic theory and use in goodness-of-fit procedures have been studied extensively \citep{BGM99,BGU05}, see \citet{PZ19} for a comprehensive review. However, the $p$-Wasserstein distance is finite only on distributions with finite $p$-th moments. When $p>1$, it therefore does not define a finite metric on the full class of distributions with finite first moments. Its empirical version also applies the $L^p$ norm directly to pointwise quantile differences, so a large discrepancy over a narrow range of upper quantiles can have a substantial effect on the empirical distance.

These limitations are especially relevant for heavy-tailed data, such as financial returns or insurance losses, where extreme observations are frequent and higher-moment assumptions are often unrealistic. This motivates a quantile-based statistic that is well defined under first-order integrability and summarizes differences cumulatively over quantile levels. To this end, we construct a new metric by repeatedly integrating the quantile functions before taking their $L^p$ distance.
We denote the resulting family of metrics by $\Delta_{n,p}$, where $n$ indexes the integrated-quantile order and $p$ is the order of the norm. When $n=1$, $\Delta_{1,p}$ is the classical $p$-Wasserstein distance; when $n\geq2$, the quantile function is integrated $n-1$ times, and $\Delta_{n,p}$ is finite under only first-order integrability. This construction is related to concepts in the stochastic-order literature: integrated distribution functions are used in stop-loss distances and higher-order stochastic dominance \citep{RR90,Rachev91,Rol76,Fish80b,SS07}, while integrated quantile functions appear in inverse stochastic dominance and related orders \citep{MS89,WY98,DR06,DC10}. These works primarily use integrated functions to formulate stochastic orders or compare risks, rather than to construct $L^p$ probability metrics under minimal moment assumptions. For $n\geq2$, the $L^p$ distance between repeatedly integrated quantile functions is finite under only finite first moments. This contrasts with the $p$-Wasserstein distance, which requires finite $p$-th moments, and is useful for two-sample testing with heavy-tailed data. We develop these integrated quantile metrics as a basis for two-sample inference and establish their metric and asymptotic properties under weak moment assumptions. 

The main contributions are as follows. First, we prove that each $\Delta_{n,p}$ is a genuine probability metric and that convergence in this metric implies weak convergence. As a byproduct, we derive a dual representation and solve an associated worst-case problem, which may be of independent interest in risk management and distributionally robust optimization. Second, for $n\ge2$, we establish strong consistency of the plug-in statistic under finite first moments, obtain its limiting distributions under the null and fixed alternatives, and provide a permutation calibration with finite-sample valid $p$-values. We further derive an asymptotic power function under local alternatives. Third, simulations show that the proposed tests achieve higher power than the Wasserstein, energy, and MMD tests under persistent same-direction quantile alternatives, while maintaining accurate size. A real data application illustrates their reduced sensitivity to extreme upper-tail observations.

The remainder of the paper is organized as follows. Section~\ref{sec-metrics} defines the integrated quantile metrics and studies their metric, convergence, and structural properties. Section~\ref{sec-inference} formulates the two-sample test and develops the asymptotic theory and permutation calibration. Section~\ref{sec-simulations} presents the simulation results. Section~\ref{sec-jtpa} gives the real data application. Section~\ref{sec-conclusion} concludes the paper. Proofs and supporting technical results are collected in the online Supplementary Material.

\textit{Notation.} We use $\Prob$, $\E$, and $\Var$ for probability, expectation, and variance, respectively; subscripts identify the relevant probability law. For a random variable $X$ with distribution function $F_X$, let $F_X^{-1}(u):=\inf\{x\in\R:F_X(x)\geq u\}$, $u\in(0,1)$, denote its left-continuous quantile function. For $1\leq r<\infty$, $X\in L^r$ means $\E[|X|^r]<\infty$, while $X\in L^\infty$ means that $X$ is essentially bounded. The symbol $\1\{\cdot\}$ denotes an indicator function, $x_+:=\max\{x,0\}$ denotes the positive part of $x$, and $\sgn(x)$ denotes its sign. For measurable functions $f$ and $g$ on $(0,1)$, let $\langle f,g\rangle:=\int_0^1 f(u)g(u)\d u$ whenever the integral is finite. We write $\|\cdot\|_p$ for the norm on $L^p(0,1)$, using the essential-supremum norm when $p=\infty$. The differential in an integral is denoted by the upright symbol $\d$. Equality and convergence in distribution are denoted by $\eqd$ and $\dto$, while $\pto$ and $\asto$ denote convergence in probability and almost sure convergence. All distributions are defined on $\R$, and the moment conditions required for each metric and result are stated explicitly.

\section{Repeatedly integrated quantile metrics}\label{sec-metrics}

This section defines the repeatedly integrated quantile metrics and establishes their main theoretical results. Section~\ref{sec2.1} proves that these distances are probability metrics and derives their basic properties. Section~\ref{sec2.2} studies convergence in these metrics, and Section~\ref{sec2.3} presents a dual representation and a related worst-case problem.

\subsection{Probability metric and basic properties}\label{sec2.1}

For an integrable random variable $X$, set $F_X^{[-1]}:=F_X^{-1}$. For $n\geq2$, define its $n$-th quantile function by repeated integration,
\begin{equation*}
\begin{aligned}
F_X^{[-n]}(u):=\int_0^u F_X^{[-n+1]}(v)\d v=\frac{1}{(n-2)!}\int_0^u F_X^{-1}(v)(u-v)^{n-2}\d v, \quad u\in[0,1].
\end{aligned}
\end{equation*}
The second equality follows from the standard identity for iterated integrals; see Proposition~2.4 of \citet{GZH26}. First-order integrability ensures that $F_X^{[-n]}$ is finite on $[0,1]$, as can be readily verified by induction on $n$ using the bound
\begin{equation*}
\left|F_X^{[-n]}(u)\right| \leq \frac{u^{n-2}}{(n-2)!} \E|X|, ~~ 0 < u < 1,\ n \geq 2.
\end{equation*}
For later use, define the linear operator $A^n:L^1(0,1)\to C^{n-2}[0,1]$ by
\begin{equation}\label{eq-24}
(A^n\varphi)(u):=\frac{1}{(n-2)!}\int_0^u\varphi(v)(u-v)^{n-2}\d v.
\end{equation}
Then $F_X^{[-n]}=A^nF_X^{-1}$ for $n\geq2$. The operator is bounded from $L^1(0,1)$ into every $L^p(0,1)$, including $p=\infty$. Its differentiability identities and exact norm bounds are collected in Section S2.1 of the Supplement.
The following definition specifies the $L^p$-distance between $n$-th quantile functions.
\begin{definition}\label{def-dist}
Let $n \geq 1$ and $p \in [1,\oo]$.
When $n=1$, assume that $X,Y \in L^p$; when $n \geq 2$, assume only that $X,Y \in L^1$.
The $L^p$-distance between the corresponding $n$-th quantile functions is defined by
\begin{equation*}
\Delta_{n,p}(X,Y) := 
    \begin{cases}
        \displaystyle 
        \left( \int_0^1 \left| F_X^{[-n]}(u) - F_Y^{[-n]}(u) \right|^p \d u \right)^{1/p}, 
        & 1 \leq p < \infty, \\[10pt]
        \displaystyle 
        \esssup_{u \in (0,1)} \left| F_X^{[-n]}(u) - F_Y^{[-n]}(u) \right|, 
        & p = \infty.
    \end{cases}
\end{equation*}
\end{definition}

Definition~\ref{def-dist} shows that, when $n=1$, $\Delta_{1,p}$ coincides with the classical $p$-Wasserstein distance. We therefore write $W_p(X,Y):=\Delta_{1,p}(X,Y)$ for this distance; see Section~1.2.3 of \citet{PZ19}.

For $n\ge2$, writing $\varphi:=F_X^{-1}-F_Y^{-1}$ and using the operator $A^n$ from \eqref{eq-24}, we have $A^n\varphi=F_X^{[-n]}-F_Y^{[-n]}$, and therefore 
$$\Delta_{n,p}(X,Y)=\bigl\| A^n(F_X^{-1}-F_Y^{-1}) \bigr\|_p.$$
Repeated integration smooths the quantile difference, as $A^n(L^1(0,1))\subset C^{n-2}[0,1]\subset L^\infty(0,1)$. It follows that $\Delta_{n,p}(X,Y)$ is finite for every $p$ whenever $X,Y\in L^1$. In contrast, for $n=1$, finiteness of the $p$-Wasserstein distance requires the stronger condition $X,Y\in L^p$.

\begin{proposition}\label{prop-metric}
For each $n \geq 1$ and $p \in [1, \oo]$, the distance $\Delta_{n, p}$ is a well-defined probability metric on its domain.
That is, for all random variables $X, Y, Z$ in the respective domains, the following hold:
\begin{enumerate}[{\rm (i)}]
\item Non-negativity and identity of indiscernibles: $\Delta_{n, p}(X, Y) \geq 0$ and $\Delta_{n, p}(X, Y) = 0 \iff X \eqd Y$.

\item Symmetry: $\Delta_{n, p}(X, Y) = \Delta_{n, p}(Y, X)$.

\item Triangle inequality: $\Delta_{n, p}(X, Z) \leq \Delta_{n, p}(X, Y) + \Delta_{n, p}(Y, Z)$.
\end{enumerate}
\end{proposition}

Proposition~\ref{prop-metric} shows that $\Delta_{n,p}$ is a genuine probability metric. In particular, the identity-of-indiscernibles property implies that $\Delta_{n,p}(X,Y)=0$ if and only if $X\eqd Y$. Thus, the repeated integration used to construct $F_X^{[-n]}$ does not make two distinct distributions indistinguishable. Consequently, $\Delta_{n,p}$ can be used as a basis for testing equality of distributions.
The boundedness of $A^n$ also yields, for $n\geq2$, $p\in[1,\oo]$, and $X,Y\in L^1$,
\begin{equation*}
\Delta_{n,p}(X,Y)\leq\frac{1}{(n-2)!}\Delta_{1,1}(X,Y)
=\frac{1}{(n-2)!}W_1(X,Y).
\end{equation*}
This inequality quantifies the smoothing effect of repeated integration and gives another direct explanation for finiteness under first-order integrability. Further monotonicity and transformation properties, together with sharper comparisons and counterexamples that clarify their scope, are given in Section S1 of the Supplement.

\subsection{Convergence properties}\label{sec2.2}

For a probability metric used in distributional inference, it is important to know whether convergence in the metric guarantees convergence of the underlying distributions. For $n=1$, convergence in the $p$-Wasserstein distance $W_p$, with $p<\oo$, is equivalent to weak convergence together with convergence of the $p$-th absolute moments; see Theorem 5.11 in \citet{Sant15}. For $n\geq2$, the bound in Section~\ref{sec2.1} shows that convergence in $W_1$ implies convergence in $\Delta_{n,p}$. The converse need not hold because repeated integration smooths differences between quantile functions. The following result shows that convergence in $\Delta_{n,p}$ is nevertheless sufficient for weak convergence.

\begin{proposition}\label{prop-weak}
Let $n \geq 2$ and $p \in [1,\oo]$. Let $\{X_m\}_{m\geq1}\subset L^1$ and $X\in L^1$.
Then 
$\Delta_{n, p}(X_m, X) \to 0 \Longrightarrow X_m \dto X.$
\end{proposition}

Proposition~\ref{prop-weak} ensures that a sequence cannot approach a target distribution in $\Delta_{n,p}$ without also converging weakly to that distribution. This property is useful for distributional estimation: convergence of an estimator in $\Delta_{n,p}$ implies weak consistency of the estimated distribution. Together with the bound in Section~\ref{sec2.1}, it gives the implication
\begin{equation*}
W_1\text{-convergence}\quad\Longrightarrow\quad
\Delta_{n,p}\text{-convergence}\quad\Longrightarrow\quad
\text{weak convergence}.
\end{equation*}
For finite $p$, and also for $p=\infty$ when $n\geq3$, convergence in $\Delta_{n,p}$ need not imply convergence in $W_1$. Weak convergence alone does not imply convergence in $\Delta_{n,p}$ for any $n\geq2$ and $p\in[1,\oo]$. The corresponding examples are given in Section S1 of the Supplement.

\subsection{Dual representation and a worst-case problem}\label{sec2.3}

In addition to the metric and convergence properties established above, this section yields two structural results of independent interest in stochastic-order analysis, risk measurement, and distributionally robust estimation and optimization.

\begin{theorem}[Dual representation]\label{thm-dual}
Let $n\geq2$, $p\in[1,\oo]$, and $X,Y\in L^1$. Let $q$ be the conjugate exponent satisfying $1/p+1/q=1$, and let $U\sim\operatorname{Unif}(0,1)$. Suppose that
$\E[F_X^{-1}(U)U^j]=\E[F_Y^{-1}(U)U^j]$ for $j=0,\ldots,n-2$.
Then 
$$\Delta_{n,p}(X,Y)
=\sup_{\phi\in\Phi}
\big|\int_0^1\{F_X^{-1}(u)-F_Y^{-1}(u)\}\phi(u)\d u\big|,$$
where
$\Phi:=\big\{\phi\in C^{n-2}[0,1]:
\phi^{(n-1)}\in L^q(0,1),\ 
\|\phi^{(n-1)}\|_q\leq1\big\},$
and $\phi^{(n-1)}$ is understood in the weak sense.
\end{theorem}

Theorem~\ref{thm-dual} expresses $\Delta_{n,p}$ through smooth weighted quantile contrasts. The moment restrictions remove the unconstrained polynomial components of the test functions; their necessity is established in Section S1.1 of the Supplement.

\begin{theorem}[Worst-case linear functional]\label{thm-worst}
Let $n\geq2$, $p\in[1,\oo)$, $X\in L^1$, and $\delta\geq0$. Let $q$ be the conjugate exponent of $p$, and let $\gamma\in C^{n-1}[0,1]$ satisfy
$\gamma^{(k)}(1)=0$ for $k=0,\ldots,n-2$.
Then
\begin{equation*}
\sup_{\substack{\varphi\in L^1(0,1)\\ \|A^n\varphi\|_p\leq\delta}}
\int_0^1\{F_X^{-1}(u)+\varphi(u)\}\gamma(u)\d u
=\int_0^1F_X^{-1}(u)\gamma(u)\d u
+\delta\|\gamma^{(n-1)}\|_q.
\end{equation*}
\end{theorem}

Theorem~\ref{thm-worst} gives the support function of an integrated-quantile perturbation ball in closed form. The perturbation is unrestricted; if $F_X^{-1}+\varphi$ is required to remain a valid quantile function, the right-hand side remains an upper bound. Proofs of both theorems are given in Section S3 of the Supplement.

\section{Two-sample testing and asymptotic theory}\label{sec-inference}

Let $X_1,\ldots,X_N\overset{\mathrm{i.i.d.}}{\sim}F_X$ and $Y_1,\ldots,Y_M\overset{\mathrm{i.i.d.}}{\sim}F_Y$ be two independent samples. We consider
\begin{equation*}
H_0:F_X=F_Y
\qquad\text{against}\qquad
H_1:F_X\neq F_Y.
\end{equation*}
By Proposition~\ref{prop-metric}, for every $n\geq2$ and $p\in[1,\oo]$, the null hypothesis holds if and only if $\Delta_{n,p}(X,Y)=0$. The case $n=1$ is the empirical $p$-Wasserstein distance, whose asymptotic theory has been studied elsewhere; see \citet{PZ19} and the references therein. We therefore restrict attention to $n\geq2$. Unless otherwise stated, $p\in[1,\oo]$. The centered fixed-alternative limit in Theorem~\ref{thm-clt} is stated for $p<\oo$; its supremum-norm counterpart requires a separate directional analysis and is not considered here.

\subsection{Test statistic and asymptotic properties}

For the $X$-sample, define the empirical distribution function
$\widehat F_{X,N}(x)={N}^{-1}\sum_{i=1}^N \1\{X_i\leq x\}, x\in\R,$
and the empirical quantile function
$\widehat F_{X, N}^{-1}(u) = \inf\{x \in \R:\ \widehat F_{X, N}(x) \geq u \}, u \in (0, 1].$
The corresponding quantities $\widehat F_{Y,M}$ and $\widehat F_{Y,M}^{-1}$ for the $Y$-sample are defined analogously.
Definition~\ref{def-dist} and the operator representation \eqref{eq-24} lead to the plug-in estimator
\begin{equation*}
\widehat \Delta_{n, p}^{(N, M)} := \left\| \widehat F_{X, N}^{[-n]} - \widehat F_{Y, M}^{[-n]} \right\|_p = \left\| A^n \(\widehat F_{X, N}^{-1} - \widehat F_{Y, M}^{-1}\) \right\|_p,
\end{equation*}
where $\widehat F_{X,N}^{[-n]}:=A^n\widehat F_{X,N}^{-1}$ and $\widehat F_{Y,M}^{[-n]}:=A^n\widehat F_{Y,M}^{-1}$ are the empirical repeatedly integrated quantile functions.

The next theorem shows that the plug-in estimator is strongly consistent under the first-order integrability condition alone.

\begin{theorem}[Strong consistency]\label{thm-cons}
If $X,Y\in L^1$, then, for every $n\geq2$ and $p\in[1,\oo]$,
$\widehat \Delta_{n, p}^{(N, M)} \asto \Delta_{n, p}(X, Y)$ as $N, M \to \oo$.
\end{theorem}

Theorem~\ref{thm-cons} establishes consistency without any smoothness or density assumptions. To derive the limiting distributions, we impose the following sufficient conditions on a generic distribution $F$.

\begin{condition}[Interior regularity and integrability]\label{cond-C1}
Let $F$ be a twice differentiable distribution function with density $f:=F'$, positive on the interior of its support $(a_F,b_F):=\{x\in\R:0<F(x)<1\}$. Assume that
\begin{equation*}
\int_0^1 \frac{t(1-t)}{f^2(F^{-1}(t))}\d t<\infty
\quad\text{and}\quad
r_F:=\sup_{0<t<1}\frac{t(1-t)|f'(F^{-1}(t))|}{f^2(F^{-1}(t))}<\infty.
\end{equation*}
\end{condition}

\begin{condition}[Endpoint regularity]\label{cond-C2}
The endpoints of the support satisfy
\begin{equation*}
a_F> -\infty
\quad\text{or}\quad
\liminf_{x\downarrow0}\frac{|f'(F^{-1}(x))|x}{f^2(F^{-1}(x))}>0,
\end{equation*}
and
\begin{equation*}
b_F<\infty
\quad\text{or}\quad
\liminf_{x\downarrow0}\frac{|f'(F^{-1}(1-x))|x}{f^2(F^{-1}(1-x))}>0.
\end{equation*}
\end{condition}

Let $B$ be a standard Brownian bridge and define
\begin{equation}\label{eq-bf}
{\mathbb B}_F(u):=\frac{B(u)}{f\{F^{-1}(u)\}},\qquad u\in(0,1).
\end{equation}
Under Conditions~\ref{cond-C1}--\ref{cond-C2}, Theorem~4.6(i) of \citet{BGU05} yields weak convergence of the empirical quantile process to $\mathbb B_F$ in $L^2(0,1)$. Lemma S2.1 in the Supplement records this result and shows that applying $A^n$ transfers the convergence to the repeatedly integrated quantile process in every $L^p(0,1)$, $p\in[1,\oo]$. This lemma is the common technical input for Theorems~\ref{thm-null} and \ref{thm-clt}; Theorem~\ref{thm-local} assumes the corresponding triangular-array quantile-process convergence. Section S4 of the Supplement verifies Conditions~\ref{cond-C1}--\ref{cond-C2} for the exponential-power and Weibull models used in the simulations.

Under $H_0$, the population quantile difference is zero, so the limiting distribution of the scaled empirical distance is a Gaussian-process norm and is generally non-normal.

\begin{theorem}[Null asymptotic distribution]\label{thm-null}
Let $F$ be a twice differentiable distribution function satisfying Conditions~\ref{cond-C1}--\ref{cond-C2}. Suppose $N,M\rightarrow \oo$, $N/(N+M)\to\lambda\in(0,1)$ and $X\eqd Y$ with common distribution $F$. For every $n\geq2$ and $p\in[1,\oo]$,
\begin{equation}\label{eq-null-limit}
r_{N,M}\widehat\Delta_{n,p}^{(N,M)}
\dto
L_{n,p,F}:=\left\|\mathbb G_{n,F}\right\|_p,
\qquad
r_{N,M}:=\sqrt{\frac{NM}{N+M}},
\end{equation}
where
$\mathbb G_{n,F}:=\sqrt{1-\lambda}\,A^n\mathbb B_F^{(1)}-\sqrt{\lambda}\,A^n\mathbb B_F^{(2)}$ and $\mathbb B_F^{(1)},\mathbb B_F^{(2)}$ are independent copies of \eqref{eq-bf}. The limit is the norm of a centered Gaussian element and is generally non-normal.
\end{theorem}

Under a fixed alternative, the population distance is positive. The form of the limit depends on $p$: the $L^p$-norm is differentiable at a nonzero argument when $p>1$, whereas the $L^1$-norm requires a directional derivative when the integrated quantile difference vanishes on part of $(0,1)$.

\begin{theorem}[Fixed-alternative asymptotic distribution]\label{thm-clt}
Let $n\geq2$ and $p\in[1,\oo)$. Let $F_X$ and $F_Y$ be twice differentiable distribution functions satisfying Conditions~\ref{cond-C1}--\ref{cond-C2}. Suppose $N,M\rightarrow \oo$,  $N/(N+M)\to\lambda\in(0,1)$ and $\Delta_{n,p}(X,Y)>0$. Set $\mu:=F_X^{[-n]}-F_Y^{[-n]}$ and define
\begin{equation*}
Z:=\sqrt{1-\lambda}\,A^n\mathbb B_{F_X}-\sqrt{\lambda}\,A^n\mathbb B_{F_Y},
\end{equation*}
where $\mathbb B_{F_X}$ and $\mathbb B_{F_Y}$ are independent. Then:
\begin{enumerate}[(i)]
\item If $p\in(1,\oo)$, define
$h^*(u):=|\mu(u)|^{p-1}\sgn\{\mu(u)\}, u\in[0,1].$
Then
\begin{equation*}
r_{N,M}\left\{\widehat\Delta_{n,p}^{(N,M)}-\Delta_{n,p}(X,Y)\right\}
\dto
\|\mu\|_p^{1-p}\langle h^*,Z\rangle
\sim \mathcal N(0,\sigma_{n,p}^2),
\end{equation*}
where
\begin{equation}\label{eq-var}
\sigma_{n,p}^2=\{\Delta_{n,p}(X,Y)\}^{2(1-p)}\Var[\langle h^*,Z\rangle].
\end{equation}
\item If $p=1$, let $E_0:=\{u\in(0,1):\mu(u)=0\}$. Then
\begin{equation*}
\begin{aligned}
r_{N,M}\left\{\widehat\Delta_{n,1}^{(N,M)}-\Delta_{n,1}(X,Y)\right\}
\dto \int_{E_0^c}\sgn\{\mu(u)\}Z(u)\d u+\int_{E_0}|Z(u)|\d u.
\end{aligned}
\end{equation*}
\end{enumerate}
\end{theorem}

Thus, the fixed-alternative limit is normal for $p>1$. For $p=1$, it is normal when $E_0$ has Lebesgue measure zero, but is generally non-normal otherwise. An explicit integral representation of $\sigma_{n,p}^2$ is given in Section S2.3 of the Supplement.

\subsection{Permutation test}
Fix a significance level $\alpha\in(0,1)$. Under $H_0$, the limiting distribution in Theorem~\ref{thm-null} depends on the unknown distribution $F$, so its critical value cannot be computed directly. We therefore use a permutation test, which proceeds as follows. First, compute the observed statistic $T_{n,p}^{\rm obs}:=\widehat\Delta_{n,p}^{(N,M)}$ from the original samples. Then pool the two samples and reassign the pooled observations into two groups of sizes $N$ and $M$. For a reassignment $\pi$, let $\widehat F_{X,N}^{\pi}$ and $\widehat F_{Y,M}^{\pi}$ be the empirical distribution functions of the two reassigned groups, and define
\begin{equation*}
T_{n,p}^{\pi}
:=\left\|A^n\left\{(\widehat F_{X,N}^{\pi})^{-1}-(\widehat F_{Y,M}^{\pi})^{-1}\right\}\right\|_p.
\end{equation*}
Enumerating all possible reassignments yields the exact permutation distribution of the test statistic, conditional on the pooled sample. Let $c_{n,p}^{\pi}(1-\alpha)$ denote its $(1-\alpha)$ quantile. The exact permutation test rejects $H_0$ if
\begin{equation*}
T_{n,p}^{\rm obs}>c_{n,p}^{\pi}(1-\alpha).
\end{equation*}
The scaling factor $r_{N,M}$ is invariant across reassignments and hence does not affect the test.
Under $H_0$, all observations have the same distribution, so the group labels are exchangeable. This gives the following result for the exact test.

\begin{proposition}\label{prop-perm-valid}
Under $H_0$, for every $n\geq2$ and $p\in[1,\oo]$, the exact
permutation test defined above satisfies $\Prob_{H_0}\left\{
T_{n,p}^{\rm obs}>c_{n,p}^{\pi}(1-\alpha)
\right\}\leq\alpha.$
\end{proposition}

Complete enumeration is usually infeasible. We therefore use a Monte Carlo permutation test based on $B$ independent reassignments $\pi_1,\ldots,\pi_B$, sampled uniformly from all possible reassignments. Its $p$-value is
\begin{equation*}
\widehat p_{n,p}^{\rm perm}
=\frac{1+\sum_{b=1}^B\1\{T_{n,p}^{\pi_b}\ge T_{n,p}^{\rm obs}\}}{B+1}.
\end{equation*}
We reject $H_0$ when $\widehat p_{n,p}^{\rm perm}\leq\alpha$. Under $H_0$, the observed labelling and the $B$ randomly generated labellings are exchangeable. Consequently, the plus-one rule gives a finite-sample level no greater than $\alpha$, with discreteness or ties making the test conservative. Thus, Proposition~\ref{prop-perm-valid} concerns complete enumeration, whereas the numerical studies use its finite-sample valid Monte Carlo version.

The following result concerns the exact permutation test. It is also consistent against every fixed alternative.
\begin{theorem}\label{thm-fixed-power}
Suppose $X,Y\in L^1$, $F_X\neq F_Y$, $N/(N+M)\to\lambda\in(0,1)$, $n\geq2$, and $p\in[1,\oo]$. Then
\begin{equation*}
T_{n,p}^{\rm obs}\asto\Delta_{n,p}(X,Y)>0,
\qquad
c_{n,p}^{\pi}(1-\alpha)\pto0.
\end{equation*}
Consequently, the exact permutation test has power converging to one.
\end{theorem}

\subsection{Local alternatives and asymptotic power}
We next study power against a sequence of alternatives approaching $H_0$. We use $K=1,2,\ldots$ to index this sequence of two-sample experiments. In experiment $K$, the two samples have sizes $N_K$ and $M_K$ and distributions $F_{X,K}$ and $F_{Y,K}$, respectively; their empirical distribution functions are denoted by $\widehat F_{X,K}$ and $\widehat F_{Y,K}$. All limits below are taken as $K\to\infty$. We assume
\begin{equation*}
N_K+M_K\longrightarrow\infty,
\qquad
\frac{N_K}{N_K+M_K}\longrightarrow\lambda\in(0,1).
\end{equation*}
Consequently, both $N_K$ and $M_K$ tend to infinity. Define
\begin{equation*}
r_K:=\left(\frac{N_KM_K}{N_K+M_K}\right)^{1/2},
\end{equation*}
so that $r_K\to\infty$. The distributions $F_{X,K}$ and $F_{Y,K}$ are allowed to vary with $K$, and their difference will be assumed to decrease at rate $r_K^{-1}$. Let $F$ be the reference distribution in the following conditions.

\begin{condition}[Local quantile-process convergence]\label{cond-C3}
Jointly in $L^2(0,1)$,
\begin{equation}\label{eq-local-process}
\sqrt{N_K}\left(\widehat F_{X,K}^{-1}-F_{X,K}^{-1}\right)\dto\mathbb B_F^{(1)},
\qquad
\sqrt{M_K}\left(\widehat F_{Y,K}^{-1}-F_{Y,K}^{-1}\right)\dto\mathbb B_F^{(2)},
\end{equation}
where $\mathbb B_F^{(1)}$ and $\mathbb B_F^{(2)}$ are independent copies of $\mathbb B_F$.
\end{condition}

\begin{condition}[Local alternative]\label{cond-C4}
For some $\eta\in L^2(0,1)$,
\begin{equation}\label{eq-local-drift}
r_K\left(F_{X,K}^{-1}-F_{Y,K}^{-1}\right)\longrightarrow\eta
\qquad\text{in }L^2(0,1).
\end{equation}
\end{condition}

Condition~\ref{cond-C3} controls the stochastic part of the local experiment: after centering at their row distributions, the two empirical quantile processes have the same Gaussian limits as under $F$. Condition~\ref{cond-C4} controls the deterministic part: the population quantile difference is of order $r_K^{-1}$, with limiting direction $\eta$. Thus, $A^n\eta$ represents the local signal, while $\mathbb G_{n,F}$ represents the sampling variation.

Condition~\ref{cond-C3} is a high-level condition for the triangular array. Lemma S2.2 of the Supplement gives sufficient conditions and covers the local quantile models used in Section~\ref{sec-simulations}. In particular, if
\begin{equation*}
F_{X,K}^{-1}=F^{-1},
\qquad
F_{Y,K}^{-1}=F^{-1}+\frac{c}{r_K}h,
\end{equation*}
where $h$ is Lipschitz and the perturbed function is a valid quantile function, then Conditions~\ref{cond-C3}--\ref{cond-C4} hold with $\eta=-ch$.

\begin{theorem}[Limit under local alternatives]\label{thm-local}
Suppose Conditions~\ref{cond-C1}--\ref{cond-C4} hold. Then, for every $n\geq2$ and $p\in[1,\oo]$,
\begin{equation*}
r_K\widehat\Delta_{n,p}^{(N_K,M_K)}
\dto
\left\|A^n\eta+\mathbb G_{n,F}\right\|_p.
\end{equation*}
\end{theorem}

Theorem \ref{thm-local} characterizes the limiting distribution of the scaled statistic under local alternatives.
Let $c_{n,p,F}(1-\alpha)$ be the $(1-\alpha)$ quantile of $L_{n,p,F}=\|\mathbb G_{n,F}\|_p$. The following corollary then translates this convergence into the asymptotic local power, and also covers the permutation version under an additional uniform continuity condition.

\begin{corollary}[Asymptotic local power]\label{cor-local-power}
Assume the conditions of Theorem~\ref{thm-local}. Suppose $c_{n,p,F}(1-\alpha)$ is the unique $(1-\alpha)$ quantile of $L_{n,p,F}$ and the distribution function of $\|A^n\eta+\mathbb G_{n,F}\|_p$ is continuous at $c_{n,p,F}(1-\alpha)$. Then
\begin{equation}\label{eq-local-power}
\begin{aligned}
\Prob\left[r_K\widehat\Delta_{n,p}^{(N_K,M_K)}>c_{n,p,F}(1-\alpha)\right]
\longrightarrow\beta_{n,p}(\eta),
\end{aligned}
\end{equation}
where $\beta_{n,p}(\eta):=\Prob\left[\|A^n\eta+\mathbb G_{n,F}\|_p>c_{n,p,F}(1-\alpha)\right]$.
Let $\mathcal Z_K$ denote the pooled data. If, in addition,
\begin{equation}\label{eq-perm-local}
\sup_{x\in\R}\left|\Prob_\pi[r_KT_{n,p}^{\pi}\leq x\mid\mathcal Z_K]-\Prob[L_{n,p,F}\leq x]\right|\pto0,
\end{equation}
then
\begin{equation*}
r_Kc_{n,p}^{\pi}(1-\alpha)\pto c_{n,p,F}(1-\alpha)
\end{equation*}
and the permutation test has limiting power $\beta_{n,p}(\eta)$.
\end{corollary}

Condition~\eqref{eq-perm-local} requires the conditional permutation distribution to converge uniformly to the null limit; standard sufficient results are given in Section~3.8 of \citet{VW23} and by \citet{CR13}. For $n=1$, \eqref{eq-local-power} gives the corresponding Wasserstein power whenever the required moment and quantile-process conditions hold. Because $A^n$ transforms both the local signal and the Gaussian noise, no general power ordering over $n$ follows.

\section{Numerical studies}\label{sec-simulations}

In this section, we examine the finite-sample performance of the proposed tests. We first assess whether the permutation tests maintain the nominal level under the null hypothesis. We then compare their power with that of $W_1$, $W_2$, energy distance, and Gaussian-kernel MMD under local alternatives. Finally, we consider a Pareto model with a finite first moment but an infinite second moment. All experiments use equal sample sizes, so $N=M$ and $\lambda=1/2$. Sections~S5.1--S5.3 of the Supplement report additional results on consistency, asymptotic approximations, and models outside the sufficient regularity conditions. Section~S5.4 gives the benchmark definitions and implementation details, while Sections~S5.5--S5.8 provide further power comparisons. Section~S5.9 provides a jointly calibrated adaptive version that combines the four proposed statistics. For simplicity, hats are omitted from empirical distance labels in figures and tables.

\subsection{Finite-sample performance of the tests}

We compare eight permutation tests at level $\alpha=0.05$: the six statistics indexed by $n\in\{1,2,3\}$ and $p\in\{1,2\}$, energy distance, and Gaussian-kernel MMD. The cases $n=1$ are $W_1$ and $W_2$, whereas $n=2,3$ give the proposed statistics. Section S5.4 of the Supplement defines the two additional benchmarks and their tuning rules.

\paragraph{Calibration under the null.}
The common law is $\mathrm{EP}_4(0,1)$ or $\mathrm{Wei}(4,1)$, with $N=M\in\{100,200,300,400,500\}$. Each setting uses $2000$ sample pairs and $199$ random relabelings per pair. For a statistic $T$, the Monte Carlo permutation $p$-value is
\begin{equation*}
\widehat p_T
=\frac{1+\sum_{b=1}^{199}\1\{T^{(b)}\geq T^{\mathrm{obs}}\}}{200}.
\end{equation*}
The null is rejected when $\widehat p_T\leq\alpha$. Exchangeability under $H_0$ and the correction by one give
$\Prob_{H_0}[\widehat p_T\leq\alpha]\leq\alpha.$
All tests use the same data and relabelings within each replication.

\begin{figure}[!htbp]
\centering
\includegraphics[width=\linewidth]{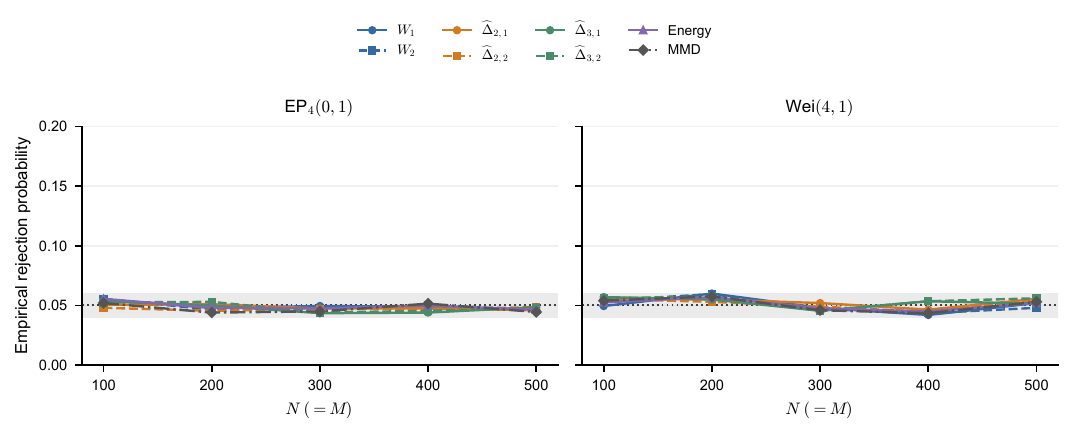}
\caption{Empirical size for $\mathrm{EP}_4(0,1)$ (left) and $\mathrm{Wei}(4,1)$ (right). The dotted line is $0.05$; the shaded band is its pointwise 95\% Monte Carlo interval.}
\label{fig-test-size}
\end{figure}

Across the 80 configurations in Figure~\ref{fig-test-size}, empirical size ranges from $0.042$ to $0.060$. All but one value lie within the pointwise Monte Carlo band, with no systematic pattern across methods or sample sizes. Results for normal, exponential, and Pareto models are given in the Supplement.

\paragraph{Power under local alternatives.}
We next consider persistent lower-quantile alternatives. Let $F=\mathrm{EP}_4(0,1)$ and $N=M=500$, so $r_{N,M}=\sqrt{250}$, and generate the second sample from
\begin{equation*}
Q_{Y,c}(u)=Q_F(u)+\frac{c}{r_{N,M}}h(u),
\qquad c\in\{0,0.4,0.8,1.2,1.6,2.0\}.
\end{equation*}
Each direction has unit $L^2(0,1)$ norm, so $c=r_{N,M}W_2(F,Q_{Y,c})$. The four directions are
\begin{equation*}
\begin{array}{ll}
\text{linear:} & h_1(u)=-\sqrt3(1-u),\\
\text{smooth:} & h_2(u)=-\sqrt{15/8}(1-u^2),\\
\text{plateau:} & h_3(u)=-\dfrac{\min\{1,(0.9-u)_+/0.2\}}{\sqrt{0.7+0.2/3}},\\[6pt]
\text{exponential decay:} & h_4(u)=-\left\{\dfrac{4}{1-e^{-4}}\right\}^{1/2}e^{-2u}.
\end{array}
\end{equation*}
The directions are nonpositive and largest near the lower endpoint, and every reported $Q_{Y,c}$ is increasing. Each setting uses $2000$ sample pairs and $199$ random relabelings. The maximum Monte Carlo standard error is $0.0112$.

\begin{figure}[H]
\centering
\includegraphics[width=\linewidth]{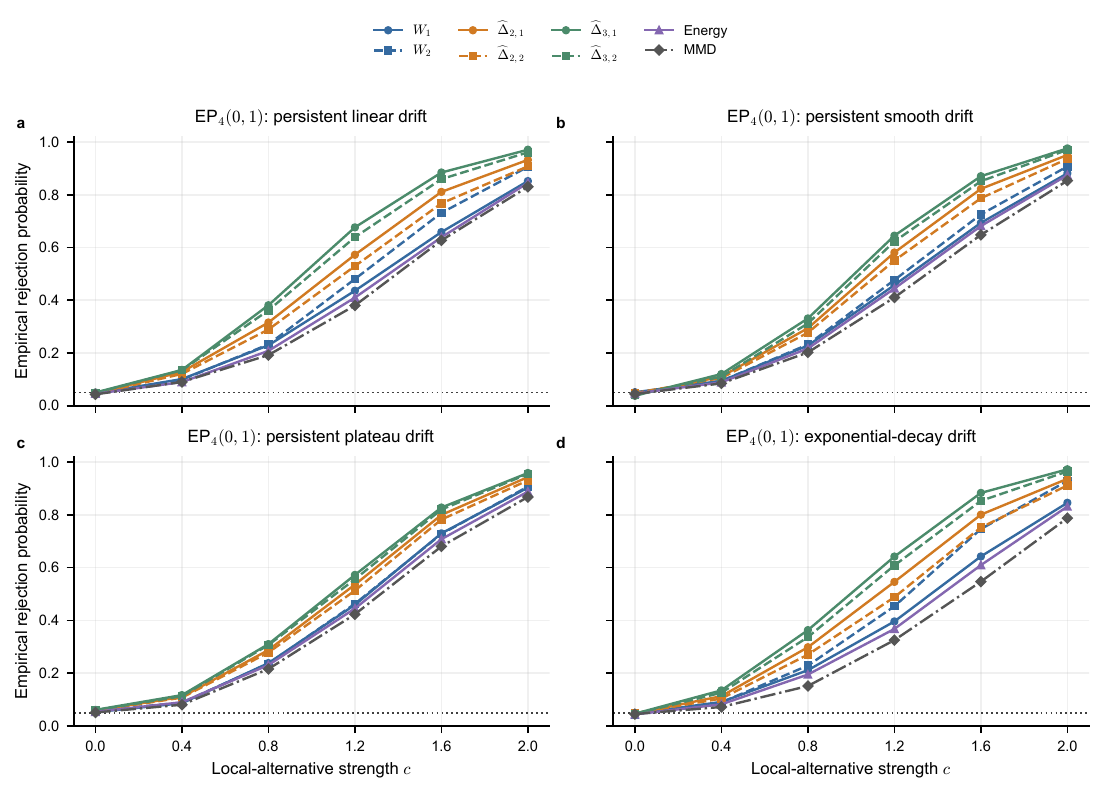}
\caption{Empirical power under the four persistent lower-quantile alternatives $h_1$--$h_4$. The dotted line is $0.05$.}
\label{fig-test-power}
\end{figure}

In Figure~\ref{fig-test-power}, every proposed test has higher empirical power than the four benchmarks for $c\in\{0.4,0.8,1.2,1.6\}$. At $c=1.2$, the strongest benchmark powers range from $0.454$ to $0.481$, whereas the proposed powers range from $0.488$ to $0.677$. The curves begin to saturate at $c=2.0$. Further base distributions and alternative shapes are reported in the Supplement.

\subsection{Simulations under weak-moment}

This section examines a setting in which the proposed metrics have finite population values under weaker moment assumptions than $W_2$. Let $F_s$ denote the Pareto$(1.5)$ distribution with quantile function
\begin{equation*}
Q_s(u)=s(1-u)^{-2/3},\qquad 0<u<1.
\end{equation*}
The distribution has a finite first moment but an infinite second moment. For any $s\neq t$,
\begin{equation*}
W_2^2(F_s,F_t)
=(s-t)^2\int_0^1(1-u)^{-4/3}\d u
=\infty,
\end{equation*}
whereas $W_1$ and the four proposed metrics considered here are finite. Although the empirical $W_2$ statistic can be computed from a finite sample, it has no finite population target in this model. Section~S5.8 of the Supplement separately studies the fixed comparison between $F_1$ and $F_{1.3}$. It reports the median and interquartile range of the sample distances as the sample size increases, together with empirical size under the Pareto null.

The main experiment focuses on two-sample testing. We draw the first sample from $F_1$ and the second from the distribution with quantile function
\begin{equation*}
Q_{Y,c}(u)
=\left(1+\frac{c}{r_{N,M}}\right)Q_1(u),
\qquad
r_{N,M}=\sqrt{\frac{NM}{N+M}}=\sqrt{250},
\end{equation*}
using $N=M=500$ and $c\in\{0,0.8,1.2,1.6\}$. Thus, $c$ controls a scale difference of order $r_{N,M}^{-1}$. For every $c>0$, the corresponding population $W_2$ distance is infinite. Each setting uses $2000$ independently generated sample pairs and $199$ random relabelings. The Pareto model is outside the sufficient quantile-process conditions used for the asymptotic approximations, so this experiment evaluates finite-sample permutation performance rather than Theorem~\ref{thm-local}.

\begin{figure}[htbp]
\centering
\includegraphics[width=0.6\linewidth]{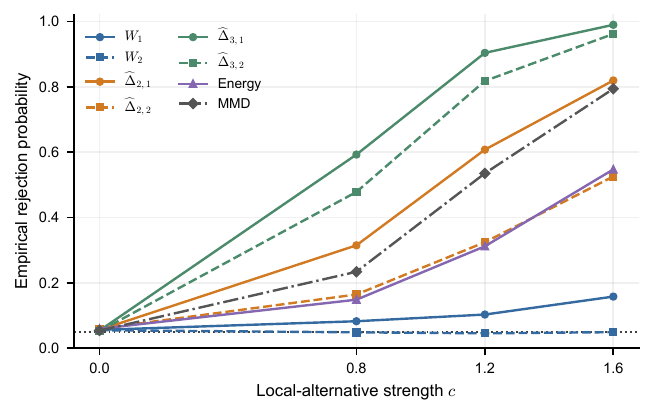}
\caption{Empirical rejection probabilities under $r_{N,M}^{-1}$-scale alternatives for the Pareto$(1.5)$ model. The dotted line is $0.05$.}
\label{fig-weak-moment-main}
\end{figure}

Figure~\ref{fig-weak-moment-main} reports the empirical rejection probabilities. At $c=1.2$, the powers of $\widehat\Delta_{2,1}$, $\widehat\Delta_{2,2}$, $\widehat\Delta_{3,1}$, and $\widehat\Delta_{3,2}$ are $0.608$, $0.324$, $0.904$, and $0.818$, respectively. The benchmark powers are $0.103$ for $W_1$, $0.045$ for $W_2$, $0.312$ for energy, and $0.535$ for MMD. Thus, the proposed tests remain applicable under first-order integrability. 

\section{Real data application}\label{sec-jtpa}

We use the National Job Training Partnership Act (JTPA) Study to examine whether assignment to an offer of JTPA services changes the distribution of subsequent earnings. We first apply the two-sample tests to pre-program earnings as a baseline-balance check and then to 30-month earnings as the primary outcome. We also compare the interpretation and upper-tail sensitivity of the proposed and benchmark statistics.

Within 16 sites of the JTPA Study, economically disadvantaged applicants were randomly assigned to an offer of employment and training services or to a control group that could not receive these services for 18 months \citep{AAI02}. Our sample contains 6102 adult women with complete 30-month earnings, including 4088 offer and 2014 control assignments. We compare groups by randomized offer, irrespective of actual participation.
The data are available from the \href{https://www.upjohn.org/data-tools/employment-research-data-center/national-jtpa-study}{Upjohn Institute's National JTPA Study archive}. We use the site-labelled adult-women extract from the \texttt{senseweight} R package and verify the outcome and assignment variables against the \href{https://economics.mit.edu/sites/default/files/publications/jtpa.raw}{MIT Economics replication file} for \cite{AAI02}. No outcome-based exclusions are made. The outcome has a point mass at zero and a long right tail, making it suitable for studying distributional differences and upper-tail sensitivity.

Let $Y_i$ denote the 30-month earnings of applicant $i$, and let $F_1$ and $F_0$ be the distributions under offer and control assignment. We test
\begin{equation*}
H_0:F_1=F_0
\qquad\text{against}\qquad
H_1:F_1\neq F_0,
\end{equation*}
using $W_1$, $W_2$, $\widehat\Delta_{2,1}$, $\widehat\Delta_{2,2}$, $\widehat\Delta_{3,1}$, $\widehat\Delta_{3,2}$, energy distance, and Gaussian-kernel MMD. Here $W_p=\widehat\Delta_{1,p}$, while $n=2,3$ gives the proposed statistics. Energy distance and MMD are defined in Section~S5.4 of the Supplement. MMD uses $512$ fixed random Fourier features \citep{RahimiRecht07} and the pooled-sample median bandwidth. The approximation differs from the exact squared MMD by 2.2\% under the observed labels.

\begin{figure}[htbp]
\centering
\includegraphics[width=\linewidth]{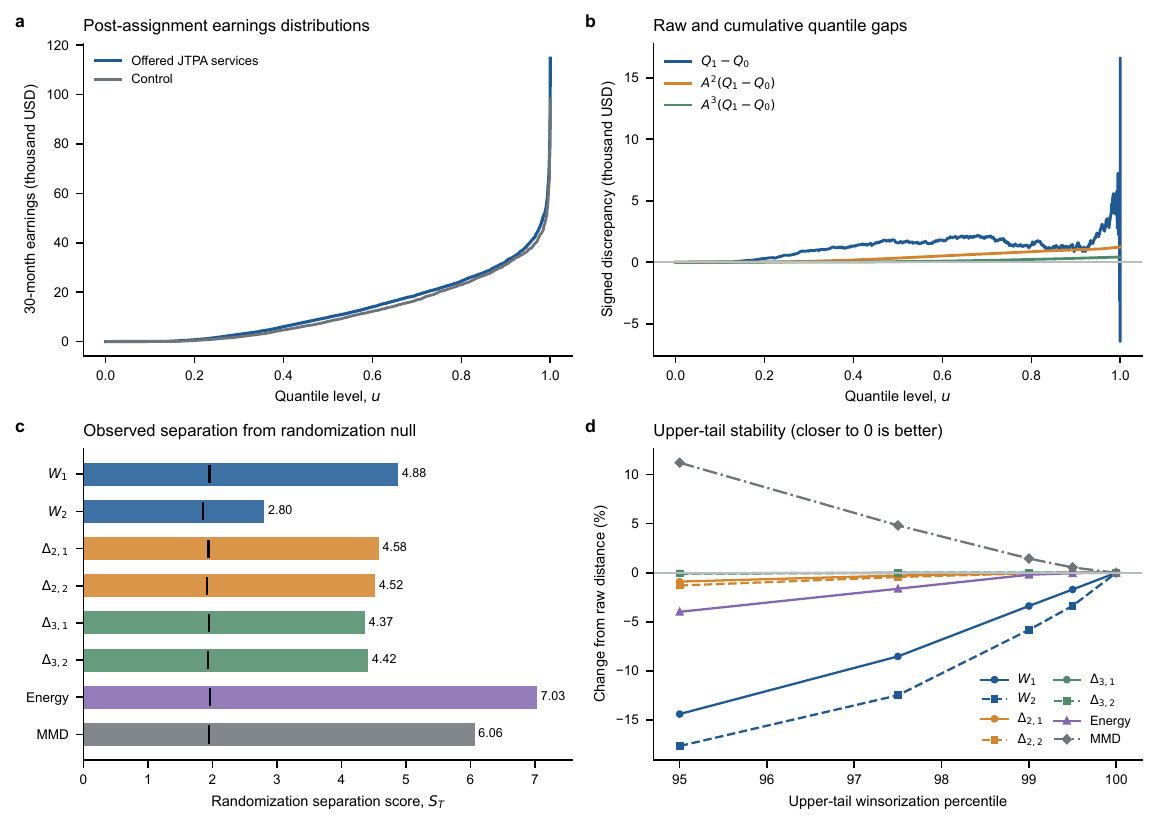}
\caption{JTPA adult-women results. \textbf{a}, Empirical quantile functions. \textbf{b}, Quantile gap $g$ and integrated profiles. \textbf{c}, Separation scores with 95\% randomization thresholds (black ticks). The scores are descriptive and do not provide a general ranking of power. \textbf{d}, Relative changes after upper-tail winsorization; zero indicates no change.}
\label{fig-jtpa}
\end{figure}

Let $\mathcal T$ denote the eight statistics, and let $T^{\rm obs}$ be the observed value of $T\in\mathcal T$. Because offers were randomized within sites, we use site-stratified randomization inference. This approach does not require the smooth-density assumptions in Section~\ref{sec-inference}, which are violated by the point mass at zero. For each of $B=9999$ draws, labels are shuffled within sites while preserving the observed group counts, and $T^{(b)}$ is recomputed. The raw randomization $p$-value is
\begin{equation*}
\widehat p_T
=\frac{1+\sum_{b=1}^{B}\1\{T^{(b)}\geq T^{\rm obs}\}}{B+1}.
\end{equation*}
Under $H_0$, the observed and shuffled labels are exchangeable within each site. The eight procedures test the same scientific null hypothesis but are reported jointly. We therefore use Holm's method to control the familywise probability of at least one false rejection across the eight reported analyses; the adjustment does not represent eight distinct scientific hypotheses \citep{Holm79}. If $\widehat p_{(1)}\leq\cdots\leq\widehat p_{(8)}$ are the ordered values, then
\begin{equation*}
\widehat p_{(i)}^{\rm Holm}
=\min\left\{1,\max_{1\leq j\leq i}(9-j)\widehat p_{(j)}\right\},
\end{equation*}
with the adjusted values returned to their original statistics. To compare statistics measured in different units, we also report
\begin{equation}\label{eq-jtpa-score}
S_T
=\frac{T^{\rm obs}-\overline T^{\rm perm}}
{s_T^{\rm perm}},
\end{equation}
where $\overline T^{\rm perm}$ and $s_T^{\rm perm}$ are the randomization mean and standard deviation. A larger $S_T$ indicates stronger separation from the statistic's own randomization distribution. The black ticks in Figure~\ref{fig-jtpa}c show the corresponding standardized 95th percentiles.

Let $\widehat Q_1$ and $\widehat Q_0$ be the empirical quantile functions and define $g=\widehat Q_1-\widehat Q_0$. We compare $g$ with
\begin{equation*}
A^2g(u)=\int_0^u g(t)\d t,
\qquad
A^3g(u)=\int_0^u (u-t)g(t)\d t.
\end{equation*}
These profiles show whether the gap accumulates over a broad quantile range. To assess upper-tail sensitivity, let $c_\tau$ be the pooled empirical $\tau$-quantile for $\tau\in\{0.95,0.975,0.99,0.995\}$ and set $Y_i^{(\tau)}=\min(Y_i,c_\tau)$. Let $T^{(\tau)}$ be the statistic computed from these winsorized outcomes under the observed labels. The MMD bandwidth and random features are held fixed. We report
\begin{equation*}
R_T^{(\tau)}
=100\left(\frac{T^{(\tau)}}{T^{\rm obs}}-1\right).
\end{equation*}
Values closer to zero indicate less upper-tail sensitivity. Figure~\ref{fig-jtpa}d gives the full path, and Table~\ref{tab-jtpa} reports the change at $\tau=0.99$.

We first apply the eight tests to pre-program earnings, using $4999$ randomizations. The raw $p$-values are $0.4042$ ($W_1$), $0.8348$ ($W_2$), $0.4742$ ($\Delta_{2,1}$), $0.5214$ ($\Delta_{2,2}$), $0.2750$ ($\Delta_{3,1}$), $0.3164$ ($\Delta_{3,2}$), $0.3774$ (energy), and $0.2706$ (MMD). All Holm-adjusted values equal one, so none of the tests rejects baseline equality.

For the primary outcome, mean earnings were USD~13,439 in the offer group and USD~12,197 in the control group. The corresponding medians were USD~9,796 and USD~8,230, while the zero-earnings shares were 12.7\% and 14.3\%. Figure~\ref{fig-jtpa}a shows that the offer-group quantiles are higher over most of the range. Although $g$ becomes erratic in the extreme upper tail, $A^2g$ and $A^3g$ remain nonnegative throughout (Figure~\ref{fig-jtpa}b). Thus, the difference accumulates over a broad quantile range.
Every score exceeds its 95\% randomization threshold in Figure~\ref{fig-jtpa}c, and all eight tests reject equality of the 30-month earnings distributions after Holm adjustment (Table~\ref{tab-jtpa}). The adjusted value is $0.012$ for $W_2$ and $0.0056$ for the other statistics. Together with the higher mean, median, and quantiles in the offer group, these results indicate a positive intention-to-treat effect of the JTPA offer on earnings.
The eight tests give the same primary conclusion, so we next compare how their statistics describe the difference. Energy and MMD have the largest separation scores, $7.03$ and $6.06$; the scores for $W_1$ and the proposed statistics range from $4.37$ to $4.88$. Each score is standardized by its own randomization distribution, whereas the raw statistics are not comparable across methods.

\begin{table}[!htbp]
\centering
\caption{Observed statistics, randomization $p$-values, and upper-tail sensitivity for 30-month earnings. The final column gives the signed percentage change after pooled 99th-percentile winsorization.}
\label{tab-jtpa}
\small
\begin{tabular}{lccccc}
\hline
Metric & \shortstack{Observed\\statistic} & \shortstack{Separation\\score $S_T$} & \shortstack{Randomization\\$p$-value} & \shortstack{Holm\\$p$-value} & \shortstack{99\% winsor\\change (\%)}\\
\hline
$W_1$ & 1247.5 & 4.88 & 0.0011 & 0.0056 & $-3.39$\\
$W_2$ & 1551.6 & 2.80 & 0.0120 & 0.0120 & $-5.82$\\
$\Delta_{2,1}$ & 423.5 & 4.58 & 0.0008 & 0.0056 & $-0.04$\\
$\Delta_{2,2}$ & 574.1 & 4.52 & 0.0007 & 0.0056 & $-0.07$\\
$\Delta_{3,1}$ & 102.6 & 4.37 & 0.0007 & 0.0056 & $-0.0006$\\
$\Delta_{3,2}$ & 160.9 & 4.42 & 0.0008 & 0.0056 & $-0.0009$\\
Energy & 74.36 & 7.03 & 0.0013 & 0.0056 & $-0.20$\\
MMD & 0.00205 & 6.06 & 0.0015 & 0.0056 & $+1.45$\\
\hline
\end{tabular}
\end{table}

Under 99th-percentile winsorization, $W_1$ and $W_2$ decrease by 3.39\% and 5.82\%, while the proposed statistics change by at most 0.07\%. Energy decreases by 0.20\%, and MMD increases by 1.45\%. At the 95th percentile, the changes reach 14.4\% for $W_1$, 17.7\% for $W_2$, 4.0\% for energy, and 11.2\% for MMD, compared with at most 1.31\% for $n=2$ and 0.10\% for $n=3$. Thus, the proposed distances are less affected by the upper-tail perturbations in this dataset.

In summary, baseline equality is not rejected, whereas all eight tests reject equality of the 30-month earnings distributions. The proposed statistics additionally show how the quantile gap accumulates and have the smallest changes under upper-tail winsorization.

\section{Conclusion}\label{sec-conclusion}

We proposed a class of two-sample tests based on repeatedly integrated quantile metrics. For $n\ge2$, these metrics are well defined under finite first moments. We established their metric properties, consistency, limiting distributions, and permutation calibration. The simulations showed accurate size and favorable power for persistent quantile differences, and the JTPA application illustrated their cumulative interpretation and reduced sensitivity to upper-tail observations. As byproducts, we also obtained a dual representation and a worst-case result.
However, the present study is restricted to univariate independent samples. \mbox{Future} work may extend the method to multivariate and high-dimensional two-sample testing and develop theory and resampling procedures for time-series and other dependent data.

\section*{Supplementary Material}
The online Supplementary Material contains additional metric properties and examples, technical proofs, regularity verification, and further simulation results.

\end{document}